\documentclass[11pt]{amsart}
\usepackage{amsmath,amsthm,amsfonts,amscd,amssymb,eucal,latexsym}

\theoremstyle{definition}
\newtheorem{theorem}{Theorem}[section]
\newtheorem{corollary}[theorem]{Corollary}
\newtheorem{lemma}[theorem]{Lemma}
\newtheorem{proposition}[theorem]{Proposition}
\newtheorem{remark}[theorem]{Remark}

\newtheorem{example}[theorem]{Example}

\newcommand{\topol}{_{\text{\rm top}}}
\newcommand{\rank}{{\rm rank}}
\newcommand{\CPA}{{\rm CPA}}
\newcommand{\rcp}{{\rm rcp}}
\newcommand{\Tr}{{\rm Tr}}

\newcommand{\Ad}{{\rm Ad}}
\newcommand{\sep}{{\rm sep}}
\newcommand{\spa}{{\rm spn}}
\newcommand{\spn}{{\rm span}}

\newcommand{\interior}{{\rm int}}
\newcommand{\Vol}{{\rm Vol}}

\allowdisplaybreaks

\begin{document}

\title[Entropy and induced dynamics]{Entropy and induced 
dynamics on state spaces}

\author{David Kerr}
\address{Mathematisches Institut, Westf\"{a}lische Wilhelms-Universit\"{a}t 
M\"{u}nster, Einsteinstr.\ 62, 48149 M\"{u}nster, Germany}
\email{kerr@math.uni-muenster.de}
\date{January 15, 2004}

\begin{abstract}
We consider the topological entropy of state space and quasi-state space
homeomorphisms induced from $C^*$-algebra automorphisms. 
Our main result asserts that, for automorphisms of separable 
exact $C^*$-algebras, zero Voiculescu-Brown entropy implies
zero topological entropy on the quasi-state space (and also more 
generally on the entire unit ball of the dual).
As an application we obtain a simple description of 
the topological Pinsker algebra in terms of local Voiculescu-Brown entropy.
\end{abstract}

\maketitle

\section{Introduction}

Entropy is a numerical invariant measuring the complexity of a dynamical 
system. For homeomorphisms of a compact metric space, topological entropy
can be defined using either open covers \cite{AKM} or separated and spanning
sets \cite{Din,Bow}. For automorphisms
of unital nuclear $C^*$-algebras, a notion of entropy based on approximation 
was introduced by Voiculescu \cite{Voi}, and this was subsequently extended
to automorphisms of exact $C^*$-algebras by Brown \cite{Br}. 
By \cite[Prop.\ 4.8]{Voi} the topological entropy
of a homeomorphism $T$ of a compact metric space $X$ coincides with the
Voiculescu-Brown entropy of the automorphism $f\mapsto f\circ T$ of
$C(X)$. In other words, the Voiculescu-Brown entropy of an automorphism
of a separable unital commutative $C^*$-algebra coincides with the topological
entropy of the induced homeomorphism on the pure state space. This paper is
principally aimed at examining the relationship 
between the Voiculescu-Brown entropy of an automorphism of a general unital
exact $C^*$-algebra and the topological entropy of the induced homeomorphism 
on the state space (or quasi-state space if we drop the requirement that
the algebra be unital) with the hope of being able to obtain information 
about one from information about the other. 

A basic problem in dynamics is to determine whether or not a given 
system has positive entropy, i.e., whether or not it is ``chaotic''.
In \cite{GW} Glasner and Weiss proved that if a homeomorphism of a 
compact metric space has zero topological entropy, then the induced
homeomorphism on the space of probability measures also has zero topological
entropy. We thus have in this case that zero Voiculescu-Brown entropy implies 
zero topological entropy on the state space. Our main result 
(Theorem~\ref{T-zero})
shows that this implication holds for automorphisms of any separable unital 
exact $C^*$-algebra, and moreover asserts that, in the general separable 
exact setting, zero Voiculescu-Brown entropy
implies zero topological entropy on the unit ball of the dual (which is
in fact equivalent to zero topological entropy on the quasi-state 
space---see Lemma~\ref{L-dual} and the paragraph preceding it).
This has the particular consequence that if an automorphism $\alpha$ of a 
$C^*$-algebra has positive topological entropy on the 
quasi-state space (for example, if $\alpha$ arises from a homeomorphism of a 
compact metric space with positive topological entropy), then any 
automorphism of a separable exact $C^*$-algebra which can be obtained from
$\alpha$ as a dynamical extension 
will have positive Voiculescu-Brown entropy. Thus in 
many cases we can obtain some information about the behaviour of 
Voiculescu-Brown entropy under taking noncommutative dynamical 
extensions, about which little seems to be known in general (it is unknown 
whether Voiculescu-Brown entropy can strictly decrease, or even become zero,
under taking dynamical extensions of a positive entropy system).

The main body of the paper is divided into three parts.
In Section~\ref{S-top} we recall the definition of topological entropy and
show that, for an automorphism of a unital $C^*$-algebra, the 
entropy of the induced homeomorphism on the state space is either zero 
or infinity. As an example we demonstrate that the shift on the full group 
$C^*$-algebra $C^* (\mathbb{F}_\infty )$ of the free group on countably 
many generators falls into the latter case. 
In Section~\ref{S-zero} we begin by recalling the definition of 
Voiculescu-Brown entropy and then proceed to the proof of our main result,
for which we develop a matrix version of an argument from \cite{GW} that 
uses results from the local theory of Banach spaces. In our case the 
key geometric fact concerns the relationship
between $n$ and $k$ given an approximately isometric embedding of $\ell^n_1$
into the space of $k\times k$ matrices with the $p=\infty$ Schatten norm. 
We round out Section~\ref{S-zero} with some applications and examples. 
In particular we show that
the shift on the reduced group $C^*$-algebra $C^*_r (\mathbb{F}_\infty )$
has zero topological entropy on the state space 
in contrast to its full group $C^*$-algebra counterpart. Finally, in 
Section~\ref{S-Pinsker} we apply our main result (or rather a local 
version that follows from the same proof) to show that the
topological Pinsker algebra from topological dynamics admits a simple
description in terms of local Voiculescu-Brown entropy. 
\medskip

\noindent{\it Acknowledgements.} 
This work was supported by the Natural
Sciences and Engineering Research Council of Canada and was carried
during stays at the University of Tokyo and the University of Rome 
``La Sapienza'' over the academic years 2001--2002 and 2002--2003,
respectively. I thank Yasuyuki Kawahigashi at the
the University of Tokyo and Claudia Pinzari at the University of Rome 
``La Sapienza'' for their generous hospitality. I 
thank Nate Brown and Hanfeng Li for helpful comments 
and Nicole Tomczak-Jaegermann for providing Lemma~\ref{L-ell1}. I also
thank one of the referees for valuable suggestions regarding 
the connection with the local theory of Banach spaces, and in particular 
for drawing my attention to the references \cite{Elt,P1,P2} and
suggesting Lemma~\ref{L-geom}, which, superceding our original lemma,
captures the geometric phenomena 
underlying Theorem~\ref{T-zero} in a general Banach space framework.

\section{Topological entropy and state space dynamics}\label{S-top}

To start with we recall the definition of topological entropy for a
homeomorphism $T$ of a compact metric space $(X,d)$ (see \cite{DGS,HK,Pet,Wal} 
for general references). For
an open cover $\mathcal{U}$ of $X$ we denote by $N(\mathcal{U})$ the
smallest cardinality of a subcover and set
$$ h_{\topol} (T,\mathcal{U}) = \lim_{n\to\infty} \frac1n \log
N(\mathcal{U}\vee T^{-1}\mathcal{U}\vee\cdots\vee T^{-(n-1)} \mathcal{U}) , $$
and we define the {\it topological entropy} of $T$ by
$$ h_{\topol}(T) = \sup_{\mathcal{U}} h_{\topol} (T,\mathcal{U}) $$
where the supremum is taken over all finite open covers $\mathcal{U}$. 
We may equivalently express the topological entropy in terms of separated
and spanning sets as follows. A set $E\subset X$ is said to be 
{\it $(n,\varepsilon )$-separated (with respect to $T$)} if for every $x,y\in 
E$ with $x\neq y$ there exists a $0\leq k \leq n-1$ such that $d(T^k x , T^k y)
> \varepsilon$, and {\it $(n,\varepsilon )$-spanning (with respect to $T$)} 
if for every $x\in X$ there is a $y\in E$ such that $d(T^k x , T^k y ) \leq 
\varepsilon$ for each $k=0, \dots , n-1$. Denoting by $\sep_n (T,\varepsilon )$
the largest cardinality of an $(n,\varepsilon )$-separated set and
by $\spa_n (T,\varepsilon )$ the smallest cardinality of an
$(n,\varepsilon )$-spanning set, we then have
$$ h_{\topol}(T) = \sup_{\varepsilon >0} \limsup_{n\to\infty}\frac1n
\log\sep_n (T,\varepsilon ) = \sup_{\varepsilon >0} \limsup_{n\to\infty}\frac1n
\log\spa_n (T,\varepsilon ) . $$

Let $A$ be a unital $C^*$-algebra. We denote by $S(A)$ the state space of $A$,
i.e., the convex set of positive linear functionals $\phi$ on $A$ with
$\phi (1) = 1$. We equip $S(A)$ with the weak$^*$ topology, under which it
is compact. Given a 
automorphism $\alpha$ of $A$ we will denote by $T_\alpha$ the
homeomorphism of $S(A)$ given by $T_\alpha (\sigma ) = \sigma\circ\alpha$
for all $\sigma\in S(A)$.

The following result is a generalization of the proposition on
p.\ 422 of \cite{Sig} (which treats homeomorphisms of compact metric spaces),
and in fact the same proof also works here. For convenience we will give a 
version of the argument in our broader context.

\begin{proposition}\label{P-di}
Let $A$ be a separable unital $C^*$-algebra and $\alpha$ an automorphism
of $A$. Then either $h_{\topol}(T_\alpha ) = 0$ or $h_{\topol}(T_\alpha ) = 
\infty$.
\end{proposition}

\begin{proof}
We define a metric $d$ on the dual $A^*$ by taking a dense sequence 
$x_1 , x_2 , x_3 , \dots$ in the unit ball of $A$ and setting
$$ d(\sigma , \omega ) = \sum_{i=1}^\infty 2^{-i} | \sigma (x_i ) - \omega
(x_i ) | . $$
The metric $d$ is compatible with the weak$^*$ topology on bounded subsets of
$A^*$. 

Now suppose that $h_{\topol}(T_\alpha ) > 0$. Then for some $\varepsilon > 0$
we have 
$$ \limsup_{n\to\infty}\frac1n \log\sep_n (T_\alpha , \varepsilon ) > 0 . $$ 
Fix an $r\in\mathbb{N}$.
Set $\theta = (1-\lambda )/(1-\lambda^r )$, and choose a $\lambda > 0$ with 
$\lambda < \min (1, 2^{-1}\theta\varepsilon )$, and set
$$ \varepsilon' = \theta\lambda^{r-1} (\varepsilon - 2\lambda\theta^{-1} ) 
> 0 . $$
For each $n\in\mathbb{N}$ let 
$E_n \subset S(A)$ be an $(n,\varepsilon )$-separated set of largest
cardinality, and consider the subset
$$ E^r_n = \bigg\{ \theta\sum_{i=1}^r \lambda^{i-1}\sigma_i : \sigma_1 ,
\dots , \sigma_r \in E_n \bigg\} $$
of $S(A)$. We will show that
$E^r_n$ is an $(n,\varepsilon' )$-separated set of cardinality $| E_n |^r$.

Indeed suppose that $(\sigma_1 , \dots , \sigma_r )$ and $(\omega_1 , \dots ,
\omega_r )$ are distinct $r$-tuples of elements in $E_n$. Then for some
$1\leq j \leq r$ we have $\sigma_i = \omega_i$ for each $i=1, \dots , j-1$
and $\sigma_j \neq \omega_j$. Since $E_n$ is $(n,\varepsilon )$-separated
there is a $0\leq k \leq n-1$ such that $d(T^k_\alpha \sigma_j ,
T^k_\alpha \omega_j ) > \varepsilon$. Observe that if $j<r$, then in
view of the definition of the metric $d$ we have
\begin{align*}
d\Big( T^k_\alpha \Big( \theta\sum_{i=j}^r \lambda^{i-1}\sigma_i \Big) , 
T^k_\alpha \big( \theta\lambda^{j-1}\sigma_j \big) \Big) &= 
d\Big( T^k_\alpha \Big( \theta\sum_{i=j+1}^r \lambda^{i-1}\sigma_i \Big) , 
0 \Big) \\
&\leq \theta\sum_{i=j+1}^r \lambda^{i-1} \\
&= \lambda^j \frac{1-\lambda^{r-j}}{1-\lambda^r} \\
&\leq \lambda^j
\end{align*}
and similarly $d\big(  T^k_\alpha \big( \theta\sum_{i=j}^r \lambda^{i-1}
\omega_i \big) , T^k_\alpha (\theta\lambda^{j-1}\omega_j )\big) 
\leq \lambda^j$. Thus when $j<r$ we have
\begin{align*}
\lefteqn{d\Big( T^k_\alpha \Big( \theta\sum_{i=1}^r \lambda^{i-1}\sigma_i 
\Big) , T^k_\alpha \Big( \theta\sum_{i=1}^r \lambda^{i-1}\omega_i \Big) \Big)}
\hspace*{2.2cm} \\ 
\hspace*{2cm} &= d\Big( T^k_\alpha \Big( \theta\sum_{i=j}^r \lambda^{i-1}
\sigma_i \Big) , T^k_\alpha \Big( \theta\sum_{i=j}^r \omega^{i-1}\sigma_i 
\Big) \Big) \\
&\geq d\big( T^k_\alpha \big( \theta\lambda^{j-1}\sigma_j \big) , T^k_\alpha
\big( \theta\lambda^{j-1}\omega_j \big) \big) -
d\Big( T^k_\alpha \Big( \theta\sum_{i=j}^r \lambda^{i-1}\sigma_i \Big) , 
T^k_\alpha \big( \theta\lambda^{j-1}\sigma_j \big) \Big) \\
&\hspace*{4cm}\ -
d\Big(  T^k_\alpha \Big( \theta\sum_{i=j}^r \lambda^{i-1}
\omega_i\Big) , T^k_\alpha \big( \theta\lambda^{j-1}\omega_j \big) \Big) \\
&\geq \theta\lambda^{j-1} d\big( T^k_\alpha \sigma_j , T^k_\alpha \omega_j 
\big) - 2\lambda^j \\
&\geq \theta\lambda^{j-1} (\varepsilon - 2\lambda\theta^{-1}) \\
&> \varepsilon' ,
\end{align*}
while in the case $j=r$ the first expression in this display is simply
$\theta\lambda^{r-1} d( T^k_\alpha \sigma_r , T^k_\alpha \omega_r )$,
which again is strictly greater than $\varepsilon'$.
It follows that $E^r_n$ is an $(n,\varepsilon' )$-separated set of 
cardinality $| E_n |^r$ as we wished to show, and so
$$ h_{\topol}(T_\alpha ) \geq \limsup_{n\to\infty}\frac1n \log\sep_n 
(T_\alpha , \varepsilon' ) \geq 
r \limsup_{n\to\infty}\frac1n \log\sep_n (T_\alpha , \varepsilon ) $$
Since $r$ was arbitrary we conclude that $h_{\topol}(T_\alpha ) = \infty$.
\end{proof}

We point out that the above argument can also be used to obtain the same 
dichotomy for the values of topological entropy among continuous maps of
state spaces induced by positive unital linear maps of separable 
operator systems or order-unit spaces, as well as among
homeomorphisms of quasi-state spaces induced by automorphisms of general
separable $C^*$-algebras (see Section~\ref{S-zero}).

As an example we will show that, for the shift on the full group 
$C^*$-algebra $C^* (\mathbb{F}_\infty )$ of the free group on countably many 
generators, the topological entropy on the state space is infinite
(Proposition~\ref{P-full}). Although 
this is an immediate consequence of the fact that the $C^*$-dynamical system 
arising from the topological $2$-shift is a $C^*$-dynamical factor of the
shift on $C^* (\mathbb{F}_\infty )$ (see the 
paragraph following Proposition~\ref{P-full}), we will give here a more
explicitly geometric proof in anticipation of the arguments in
Section~\ref{S-zero}. In addition, the following two lemmas which we will 
require are of use in other situations. For example, 
Lemma~\ref{L-dual} implies that the values of topological entropy
on the state space and quasi-state space agree (see the second paragraph
of Section~\ref{S-zero}), and its proof can also be used to show that the 
values of topological entropy on the quasi-state space and on the unit ball
of the dual agree. Lemma~\ref{L-basis}, on the other hand, 
will be of convenience in Example~\ref{E-operator}.

\begin{lemma}\label{L-dual}
Let $A$ be a separable unital $C^*$-algebra and $\alpha$ an automorphism
of $A$. Let $S_\alpha$ be the homeomorphism of the closed unit ball 
$B_1 (A^* )$ of the dual of $A$ (with the weak$^*$ topology) 
given by $S_\alpha (\sigma ) = \sigma\circ\alpha$. Then
$$ h_{\topol}(S_\alpha ) = h_{\topol}(T_\alpha ) . $$
\end{lemma}

\begin{proof}
Since $T_\alpha$ is the restriction of $S_\alpha$ to $S(A)$, we have
$h_{\topol}(S_\alpha ) \geq h_{\topol}(T_\alpha )$. Thus, in view of 
Proposition~\ref{P-di}, we need only show that $h_{\topol}(T_\alpha ) = 0$
implies $h_{\topol}(S_\alpha ) = 0$. Suppose then that
$h_{\topol}(T_\alpha ) = 0$. Let $x_1 , x_2 , x_3 , \dots$ be
a dense sequence in the unit ball of $A$ and define on $A^*$ the metric
$$ d(\sigma , \omega ) = \sum_{i=1}^\infty 2^{-i} | \sigma (x_i ) - \omega
(x_i ) | , $$
which is compatible with the weak$^*$ topology on bounded subsets of
$A^*$. Let $\varepsilon > 0$, and
pick an integer $r > \varepsilon^{-1}$. Let $E_n \subset S(A)$ be an 
$(n,\varepsilon )$-spanning set (with respect to $T_\alpha$) of smallest 
cardinality. Now if $\tau\in B_1 (A^* )$ then we can write
$$ \tau = \sigma_1 - \sigma_2 + i(\sigma_3 - \sigma_4 ) $$
where each $\sigma_j$ is a positive linear functional of norm at most $1$
\cite[Thm.\ 4.3.6 and Cor.\ 4.3.7]{KR}. Since $E_n$ is 
$(n,\varepsilon )$-spanning, for each $j=1,2,3,4$ for which $\sigma_j \neq 0$
we can find an $\omega_j \in E_n$ such that
$$ d\big( T^k_\alpha \omega_j , T^k_\alpha (\sigma_j / \| \sigma_j \|)\big)
\leq\varepsilon $$
for every $k=0, \dots , n-1$, and we can also find an
$m_j \in\{ 0,1,\dots , r\}$ such that 
$\big| m_j /r - \| \sigma_j \| \big| \leq\varepsilon$, so that, 
for every $k=0, \dots , n-1$, 
\begin{align*}
d\big( S^k_\alpha ((m_j /r)\omega_j ) , S^k_\alpha \sigma_j \big) 
&\leq d\big( S^k_\alpha ((m_j /r)\omega_j ) , S^k_\alpha (\| \sigma_j \|
\omega_j ) \big) \\
&\hspace*{3cm}\ + d\big( S^k_\alpha (\| \sigma_j \| \omega_j ) ,
S^k_\alpha \sigma_j \big) \\
&\leq \big| m_j /r - \| \sigma_j \| \big| + \| \sigma_j \| d\big(  
T^k_\alpha \omega_j , T^k_\alpha (\sigma_j / \| \sigma_j \|)\big) \\
&\leq \varepsilon + \| \sigma_j \| \varepsilon \\
&\leq 2\varepsilon . 
\end{align*}
For any $j=1,2,3,4$ for which $\sigma_j = 0$, let $\omega_j$ be any
bounded linear functional on $A$ and set $m_j = 0$.
Set $\tau' = (m_1 \omega_1 - m_2 \omega_2 + im_3 \omega_3 - im_4 
\omega_4 )/r$. Then $\| \tau' \| \leq 4$, and if $\tilde{S}_\alpha$ denotes
the homeomorphism of the closed ball $B_4 (A^* ) = \{ \sigma\in A^* :
\| \sigma \| \leq 4 \}$ given by $\tilde{S}_\alpha (\sigma ) = 
\sigma\circ\alpha$, then in view of the definition of the metric $d$ we have
$$ d\big( \tilde{S}^k_\alpha \tau' , \tilde{S}^k_\alpha \tau \big) \leq 
\sum_{i=1}^4 d\big( S^k_\alpha ((m_i /r)\omega_i ) , S^k_\alpha \sigma_i \big)
\leq 8\varepsilon + 8\varepsilon = 16\varepsilon $$
for all $k=0, \dots , n-1$. Let $F_n$ be the subset of $B_4 (A^* )$ 
consisting of all $\tau'$ which arise in the 
above way with respect to some $\tau\in B_1 (A^* )$. We have thus shown
that $F_n$ is $(n,16\varepsilon )$-spanning for $B_1 (A^* )$ with respect to
$\tilde{S}_\alpha$ (in the obvious sense which relativizes the definition
of an $(n,\varepsilon )$-spanning set to a subset of the space---see 
Definition 14.14 of \cite{DGS}), and hence that the smallest cardinality 
$\spa_n (\tilde{S}_\alpha , 16\varepsilon , B_1 (A^* ))$ of an 
$(n,16\varepsilon )$-spanning set for $B_1 (A^* )$ with respect to 
$\tilde{S}_\alpha$ is bounded above by $| F_n |$, which is in turn bounded
above by $| E_n |^{4(r+1)}$.
Since $r$ does not depend on $n$ we therefore obtain
$$ \limsup_{n\to\infty}\frac1n \log\spa_n (\tilde{S}_\alpha , 16\varepsilon , 
B_1 (A^* )) \leq 4(r+1)\limsup_{n\to\infty}\frac1n \log | E_n | = 0 . $$
Now Proposition 14.15 of \cite{DGS} shows that the supremum of the first 
expression in the above display over all $\varepsilon > 0$ is in fact 
equal to the topological entropy of $S_\alpha = 
\tilde{S}_\alpha \big| {}_{B_1 (A^* )}$, and so we conclude that
$h_{\topol}(S_\alpha ) = 0$, as desired.
\end{proof}

\begin{lemma}\label{L-basis}
Let $A$ be a unital $C^*$-algebra and $\alpha$ an automorphism of $A$. 
Suppose that there exists an $x\in A$ and a $K\geq 1$ such that, for each
$n\in\mathbb{N}$, the linear map $\Gamma_n : \ell^n_1 \to\spn 
\{ x , \alpha (x) , \dots , \alpha^{n-1}(x) \}$ which sends
the $i$th standard basis element of $\ell^n_1$ to $\alpha^{i-1} (x)$ for 
each $i=1, \dots , n$ is an isomorphism 
whose inverse has norm at most $K$. Then $h_{\topol}(T_\alpha ) = \infty$.
\end{lemma}

\begin{proof}
Let $x\in A$, $K\geq 1$, and $\Gamma_n$ for $n\in\mathbb{N}$ be as in the 
hypotheses. For each $n\in\mathbb{N}$ we denote by
$\Lambda_n$ the collection of functions from $\{ 0 , \dots , n-1 \}$ to
$\{ -1 , 1 \}$. For each $f\in\Lambda_n$ we define the linear functional
$\sigma_f$ on $\spn (x, \alpha (x) , \dots , \alpha^{n-1}(x) )$ by
specifying
$$ \sigma_f \big( \alpha^k (x)\big) = K^{-1} f(k) $$
for each $k=0, \dots , n-1$. Since $\Gamma^{-1}_n$ has norm at most $K$,
it follows that $\sigma_f$ has norm at most one.
By the Hahn-Banach theorem we can extend $\sigma_f$ to an element 
$\sigma_f'$ in the dual $A^*$ of norm at most one. Let $\mathcal{U}$
be the open cover of the closed unit ball $B_1 (A^* )$ consisting of
the two open sets
\begin{gather*}
\{ \sigma\in B_1 (A^* ) : | \sigma (x) - K^{-1} | < 2K^{-1} \} , \\
\{ \sigma\in B_1 (A^* ) : | \sigma (x) - K^{-1} | > K^{-1} \} .
\end{gather*}
Let $S_\alpha : B_1 (A^* ) \to B_1 (A^* )$ be the 
homeomorphism given by $S_\alpha (\sigma ) = \sigma\circ\alpha$. For every 
$n\in\mathbb{N}$, each element of $\mathcal{U}\vee S^{-1}_\alpha\mathcal{U}
\vee\cdots\vee S^{-(n-1)}_\alpha\mathcal{U}$ contains precisely one linear
functional of the form $\sigma_f$ for $f\in\Lambda_n$, from which it follows
that 
$$ N(\mathcal{U}\vee S^{-1}_\alpha\mathcal{U}\vee\cdots\vee S^{-(n-1)}_\alpha
\mathcal{U}) = 2^n . $$
Consequently
$$ h_{\topol}(S_\alpha ) \geq\lim_{n\to\infty}\frac1n \log N(\mathcal{U}\vee 
S^{-1}_\alpha\mathcal{U}\vee\cdots\vee S^{-(n-1)}_\alpha \mathcal{U}) = 
\log 2 , $$
and so $h_{\topol}(T_\alpha ) = \infty$ by Lemma~\ref{L-dual} and
Proposition~\ref{P-di}.
\end{proof}

Let $\{ u_i \}_{i\in\mathbb{Z}}$ be 
the set of canonical unitaries associated to the generators 
in the full group $C^*$-algebra $C^* (\mathbb{F}_\infty )$. The shift 
$\alpha$ on $C^* (\mathbb{F}_\infty )$ is the automorphism defined by 
specifying $\alpha (u_i ) = u_{i+1}$ for all $i\in\mathbb{Z}$.

\begin{proposition}\label{P-full}
For the shift $\alpha$ on $C^* (\mathbb{F}_\infty )$ we have 
$$ h_{\topol}(T_\alpha ) = \infty . $$
\end{proposition}

\begin{proof}
For each $n\in\mathbb{N}$, the linear map from $\spn (u_0 , \dots , 
u_{n-1} )$ to $\ell^n_1$ which sends $u_i$ to the $i$th
standard basis element of $\ell^n_1$ is an isometry, as pointed out
in \cite[Sect.\ 8]{Pis}. This follows from
the observation that, for any scalars $c_0 , \dots , c_{n-1}$, 
$$ \left\| \sum_{i=0}^{n-1} c_i u_i \right\| = 
\sup \left\| \sum_{i=0}^{n-1} c_i v_i \right\| $$
where the supremum is taken over all unitaries $v_0 , \dots , v_{n-1}\in
\mathcal{B}(\ell_2 )$. We can thus appeal to 
Lemma~\ref{L-basis} to obtain the result.
\end{proof}

Proposition~\ref{P-full} can also be established by observing that if
$T$ is the left shift on $X = \{ -1,1 \}^{\mathbb{Z}}$ and $\beta$ is the
automorphism of $C(X)$ defined by $\beta (f) = f\circ T$ for all 
$f\in C(X)$, then $\beta$ is a $C^*$-dynamical factor of 
$\alpha$. Indeed if we consider for each $n\in\mathbb{N}$ the function 
$g_n \in C(X)$ given by
$$ g_n ((a_k )_{k\in\mathbb{Z}}) = a_n $$
for all $(a_k )_{k\in\mathbb{Z}} \in X$, then we can define a 
$^*$-homomorphism $\gamma : C^* (\mathbb{F}_\infty ) \to C(X)$ by
specifying $\gamma (u_n ) = g_n$ for each $n\in\mathbb{N}$. By the
Stone-Weierstrass theorem $\gamma$ is surjective, and evidently 
$\gamma\circ\alpha = \beta\circ\gamma$. It follows that $T_\alpha$ contains
a subsystem conjugate to $T_\beta$ and hence has infinite entropy by
Proposition~\ref{P-di}.

In contrast we will show in Proposition~\ref{P-red} that, for the 
corresponding shift on the reduced group $C^*$-algebra $C^*_r 
(\mathbb{F}_\infty )$, the topological entropy on the state space is zero.

\section{Zero Voiculescu-Brown entropy implies zero topological 
entropy on the unit ball of the dual}\label{S-zero}

We begin this section by recalling the definition of 
Voiculescu-Brown entropy, which is an extension to exact $C^*$-algebras 
\cite{Br} of the approximation-based entropy for automorphisms of unital 
nuclear $C^*$-algebras introduced in \cite{Voi}.
Let $A$ be an exact $C^*$-algebra, and let $\pi : A\to\mathcal{B}
(\mathcal{H})$ be a faithful representation. For a finite set 
$\Omega\subset A$ and $\delta > 0$ we denote by 
$\CPA (\pi , \Omega , \delta )$ 
the collection of triples $(\phi , \psi , B)$ where $B$ is a 
finite-dimensional $C^*$-algebra and $\phi : A\to B$ and $\psi : B\to
\mathcal{B}(\mathcal{H})$ are contractive completely positive linear maps 
such that
$\| (\psi\circ\phi )(x) - \pi (x) \| < \delta$ for all $x\in\Omega$. This
collection is non-empty by nuclear embeddability \cite{Kir}. We 
define $\rcp (\Omega , \delta )$ to be the infimum of $\rank\, B$ over
all $(\phi , \psi , B)\in\CPA (\pi , \Omega , \delta )$, where rank refers
to the dimension of a maximal Abelian $C^*$-subalgebra. As the notation
indicates, this infimum is independent of the particular faithful
representation $\pi$, as demonstrated in the proof of 
Proposition 1.3 in \cite{Br}. For an 
automorphism $\alpha$ of $A$ we set 
\begin{align*}
ht(\alpha , \Omega , \delta ) &= \limsup_{n\to\infty}\frac1n \log\rcp 
(\Omega\cup\alpha\Omega\cup\cdots\cup\alpha^{n-1}\Omega , \delta ) , \\
ht(\alpha , \Omega ) &= \sup_{\delta > 0}ht(\alpha , \Omega , 
\delta ) , \\
ht(\alpha ) &= \sup_\Omega ht(\alpha , \Omega ) 
\end{align*}
with the last supremum taken over all finite sets $\Omega\subset A$. We refer
to $ht(\alpha )$ as the {\it Voiculescu-Brown entropy} of $\alpha$.

For any $C^*$-algebra $A$, the quasi-state space $Q(A)$ of $A$ is 
defined as the convex set of positive linear functionals $\phi$ on $A$ with
$\| \phi \| \leq 1$. 
Equipped with the weak$^*$-topology, $Q(A)$ is compact.
Given an automorphism $\alpha$ of $A$ we denote by $\tilde{T}_\alpha$
the homeomorphism $\sigma\mapsto\sigma\circ\alpha$ of $Q(A)$. 
As in the previous section, when $A$ is unital we denote by $T_\alpha$ the 
homeomorphism $\sigma\mapsto\sigma\circ\alpha$ of the state space $S(A)$,
(i.e., the restriction of $\tilde{T}_\alpha$ to $S(A)$), and for
general $A$ we denote by $S_\alpha$ the homeomorphism 
$\sigma\mapsto\sigma\circ\alpha$ of the closed unit ball $B_1 (A^* )$ of 
the dual of $A$ with the weak$^*$ topology.
Note that, since $Q(A)$ is a subset of the unit ball of 
the dual $A^*$, by Lemma~\ref{L-dual} we have $h_{\topol}(\tilde{T}_\alpha ) 
= h_{\topol}(T_\alpha )$ in the separable unital case. The argument 
in the proof of Lemma~\ref{L-dual} can also be used to show 
that $h_{\topol}(S_\alpha ) = h_{\topol}(\tilde{T}_\alpha )$ in the
general separable case. 

Before coming to the statement of our main result we establish two
lemmas. These involve problems of a typical nature
in the local theory of Banach spaces and will be proved using methods
from this theory. For the basic background we refer the reader to \cite{MS}
and \cite{BM}. Here we may take our Banach spaces to be over either 
the real or complex numbers, but
it is the complex case which is relevant for our applications. For 
$1\leq p \leq\infty$ we denote by $C^k_p$ the
Schatten $p$-class, i.e., the space of $k\times k$ matrices with norm
$\| x \|_p = \Tr (|x|^p )^{1/p}$ in the case $1\leq p < \infty$ (where
$\Tr$ is the trace taking value $1$ on minimal projections), or
the operator norm (with the matrices operating on $\ell^k_2$) 
in the case $p=\infty$.
Given isomorphic Banach spaces $X$ and $Y$ and $K\geq 1$, we say that $X$ and 
$Y$ are {\it $K$-isomorphic} if the Banach-Mazur distance
$$ d(X,Y) = \inf \{ \| \Gamma \| \| \Gamma^{-1} \| : \Gamma :X\to Y 
\text{ is an isomorphism} \} $$
is no greater than $K$. 

The following lemma and its proof were communicated to me by
Nicole Tomczak-Jaegermann.

\begin{lemma}\label{L-ell1}
Let $X$ be an $n$-dimensional subspace of $C^k_\infty$ which is 
$K$-isomorphic to $\ell^n_1$. Then 
$$ n \leq a K^2 \log k $$
where $a>0$ is a universal constant.
\end{lemma}

\begin{proof}
The idea, which is standard in the local theory of Banach spaces, is to 
compare (Rademacher) type $2$ constants. It can be seen from the
proof of Theorem 3.1(ii) in \cite{TC} that the type $2$ constant of the 
Schatten $p$-class for $2\leq p < \infty$ satisfies
$$ T_2 (C^k_p ) \leq C \sqrt{p} $$
where $C$ is a universal constant. Thus, since the Banach-Mazur distance
$d(C^k_\infty , C^k_p )$ is equal to $k^{1/p}$ for $2\leq p < \infty$
\cite[Thm.\ 45.2]{BM}, we have
$$ T_2 (C^k_\infty ) \leq d(C^k_\infty , C^k_p ) T_2 (C^k_p ) \leq 
Ck^{1/p}\sqrt{p} $$
for every $2\leq p < \infty$. Setting $p = \log k$ we obtain for sufficiently
large $k$ the bound
$$ T_2 (C^k_\infty ) \leq Ce \sqrt{\log k} , $$
and since the type $2$ constant for $\ell^n_1$ satisfies 
$T_2 (\ell^n_1 ) \geq\sqrt{n}$ (see \S 4 in \cite{BM}) it follows that
$$ \sqrt{n} \leq T_2 (\ell^n_1 ) \leq K T_2 (X) \leq K 
T_2 (C^k_\infty )\leq KCe \sqrt{\log k} , $$
yielding the assertion of the lemma.
\end{proof}

The next lemma is a matrix analogue of Proposition 2.1 of \cite{GW}.  
I thank one of the referees for suggesting its formulation as a general 
Banach space result. The main part of its proof is concerned with
symmetric convex subsets of the unit cube which contain many 
$\varepsilon$-separated points, and was first 
established by Elton
\cite{Elt} in the real case and Pajor \cite{P1,P2} in the complex case.
Our general line of argument follows \cite{GW}, with the part 
involving almost Hilbertian sections of unit balls being replaced in our case
with an appeal to Lemma~\ref{L-ell1}.
As is our convention, for a Banach space $X$ we write $B_r (X)$ to refer 
to the closed ball $\{ x\in X : \| x \| \leq r \}$.

\begin{lemma}\label{L-geom}
Given $\varepsilon > 0$ and $\lambda > 0$ there exists a $\mu > 0$ such 
that, for all $n\geq 1$, if $\phi : C_1^{r_n} \to
\ell_\infty^n$ is a contractive linear map such that 
$\phi (B_1 (C_1^{r_n}))$ contains an $\varepsilon$-separated set of
cardinality at least $e^{\lambda n}$, 
then $r_n \geq e^{\mu n}$. 
\end{lemma}

\begin{proof}
It is well known in Banach space theory that for every $\varepsilon > 0$
and $\lambda > 0$ there exist $d > 0$ and $\delta > 0$ such that the
following holds for all $n\geq 1$: if $S\subset B_1 (\ell_\infty^n )$ 
is a symmetric convex set which contains an $\varepsilon$-separated set $F$
of cardinality at least $e^{\lambda n}$ then there is a subset
$I_n \subset \{ 1,2,\dots , n \}$ with cardinality at least $dn$ such 
that 
$$ B_\delta (\ell_\infty^{I_n}) \subset \pi_n (S) $$
where $\pi_n : \ell_\infty^n \to \ell_\infty^{I_n}$ is the canonical
projection. In the real case this is implicit in the argument on p.\ 117
in \cite{Elt}, while in the complex case it follows immediately from
Th\'{e}or\`{e}me 5 of \cite{P2} (take $t=\varepsilon /2$; then, with  
$K = B_1 (\ell_\infty^n )$, the cubes $x+tK$ for $x\in F$ have 
pairwise disjoint interiors, whence 
$\Vol (S+tK) \geq t^n 2^n e^{\lambda n}$, yielding the conclusion with,
e.g., $\delta = (e^\lambda - 1)\varepsilon /4$, assuming 
$\varepsilon\leq 1/2$).
Thus taking $S = \phi (B_1 (C_1^{r_n}))$ we obtain
$$ B_\delta (\ell_\infty^{I_n}) \subset \pi_n (\phi (B_1 (C_1^{r_n}))) . $$
Hence the dual map $(\pi_n \circ\phi )^*$ from
$(\ell_\infty^{I_n})^* \cong\ell_1^{I_n}$ to 
$(C_1^{r_n})^* \cong C_\infty^{r_n}$ is an embedding of norm at most $1$
whose inverse has norm at most $1/\delta$. Since these bounds do not 
depend on $n$, by Lemma~\ref{L-ell1} there 
is a $c>0$ such that, for all $n\geq 1$, 
$$ | I_n | \leq c \log r_n $$
and hence $dn \leq c \log r_n$. Setting $\mu = d/c$ we obtain the
assertion of the lemma.
\end{proof}

\begin{theorem}\label{T-zero}
Let $A$ be a separable exact $C^*$-algebra and $\alpha$ an
automorphism of $A$. Then $ht(\alpha ) = 0$ implies $h_{\topol}(S_\alpha )
=0$, and hence also $h_{\topol}(\tilde{T}_\alpha ) = 0$ and (when $A$
is unital) $h_{\topol}(T_\alpha ) = 0$. 
\end{theorem}

\begin{proof}
Let $K$ be any compact subset of the unit ball of 
$A$ whose linear span is dense in $A$ (for example, we may take 
$K = \{ k^{-1}x_k \}_{k\in\mathbb{N}}$ where $\{ x_k \}_{k\in\mathbb{N}}$ 
is a dense sequence in the unit ball of $A$). 
The metric $d$ on $B_1 (A^* )$ defined by
$$ d (\sigma , \omega ) = \sup_{x\in K} | \sigma (x) - \omega (x) | $$
for all $\sigma , \omega\in B_1 (A^* )$ is readily seen to give rise to the 
weak$^*$ topology. Now suppose $h_{\topol}(S_\alpha ) > 0$. 
Then there exist an $\varepsilon > 0$, a $\lambda > 0$, and an infinite
set $J\subset\mathbb{N}$ such that for all $n\in J$ there is an 
$(n,4\varepsilon )$-separated set $E_n \subset B_1 (A^* )$ of cardinality at 
least $e^{\lambda n}$. By compactness there is a finite set $\Omega\subset
K$ such that, for all $\sigma , \omega\in B_1 (A^* )$,
$$ d (\sigma , \omega ) \leq \sup_{x\in\Omega} | \sigma (x) - 
\omega (x) | + \varepsilon . $$
We will show that $ht(\alpha , \Omega , \varepsilon ) > 0$.

Let $\pi : A \to\mathcal{B}(\mathcal{H})$ be any faithful representation.
For each $n\in J$, let $(\phi_n , \psi_n , B_n )$ be an element in 
$\CPA (\pi , \Omega\cup\cdots\cup\alpha^{n-1}\Omega , \varepsilon )$ 
with $B_n$ of smallest possible rank, and set $r_n = \rank\, B_n$. 
Writing $\Omega = \{ x_1 , \dots , x_m \}$ we define a map $\Gamma_n$ from 
the Schatten class $C_1^{r_n}$ to $(\ell_\infty^n )^m \cong\ell_\infty^{nm}$ 
by
$$ \Gamma_n (h) = (( \Tr (h\phi_n (\alpha^k (x_i ))))_{k=1}^{n-1}
)_{i=1}^{m} $$ 
for all $h\in C_1^{r_n}$, where $\Tr$ is the trace on $M_{r_n} (\mathbb{C})$ 
taking value $1$ on minimal projections and $B_n$ is considered as a
$C^*$-subalgebra of $M_{r_n} (\mathbb{C})$ under some fixed embedding. 
Note that $\Gamma_n$ is contractive, since
$\Omega\cup\alpha\Omega\cup\dots\cup\alpha^{n-1}\Omega$ lies in the unit 
ball of $A$, $\phi_n$ is contractive, and 
$|\Tr (hx)| \leq \Tr (|h|) \| x \|$ for all $h\in C_1^{r_n}$ and $x\in B_n$. 

For each $\sigma\in B_1 (A^* )$ we can extend $\sigma\circ\pi^{-1}$ on 
$\pi (A)$ to a contractive linear functional $\sigma'$ on 
$\mathcal{B}(\mathcal{H})$ by the Hahn-Banach theorem. 
Now if $\sigma$ and $\omega$ are distinct 
elements of $E_n$ then there is a $k$ with $0\leq k\leq n-1$ such that 
$d (T_\alpha^k \sigma , T_\alpha^k \omega ) > 4\varepsilon$. Then
$$ \sup_{x\in\Omega} \big| (\sigma\circ\alpha^k )(x) - (\omega\circ\alpha^k
)(x) \big| > 3\varepsilon , $$
and since for every $x\in\Omega$ we have
$$ \| (\psi_n \circ\phi_n )(\alpha^k (x)) - \pi (\alpha^k (x)) \| < 
\varepsilon $$
it follows by the triangle inequality that
$$ \sup_{x\in\Omega} \big| (\sigma' \circ\psi_n ) (\phi_n (\alpha^k (x)))
- (\omega' \circ\psi_n ) (\phi_n (\alpha^k (x))) \big| > \varepsilon . $$
Taking a conditional expectation $P : M_{r_n} (\mathbb{C})\to B_n$ and
isometrically identifying a linear functional on 
$M_{r_n} (\mathbb{C})$ with its density matrix in $C_1^{r_n}$, we thus have 
that the image of the subset
$\{ \sigma' \circ\psi_n \circ P : \sigma\in E_n \}$ of $B_1 (C_1^{r_n})$  
under $\Gamma_n$ is an
$\varepsilon$-separated set with cardinality at 
least $e^{\lambda n}$. Since $m$ does not depend on $n$,
by Lemma~\ref{L-geom} there is a $\mu > 0$ such that $r_n \geq e^{\mu n}$ for 
all $n\in J$, and so 
$ht(\alpha , \Omega , \varepsilon )\geq \mu > 0$, yielding the result.
\end{proof}

\begin{remark}
Notice that in the proof of Theorem~\ref{T-zero} we made no use of the
order structure. In fact Pop and Smith have shown in \cite{SP} 
that Voiculescu-Brown entropy can be alternatively defined using 
completely contractive linear maps. 
\end{remark}

The following corollary is an immediate consequence of 
Theorem~\ref{T-zero}
and the fact that topological entropy does not increase under taking
factors or restrictions to closed invariant subsets.

\begin{corollary}\label{C-pos}
Let $A$ and $B$ be separable exact $C^*$-algebras and 
$\alpha :A\to A$ and
$\beta :B\to B$ automorphisms with $h_{\topol}(S_\alpha ) > 0$.
Suppose that there exists a surjective contractive linear map
$\gamma : B\to A$ such that $\alpha\circ\gamma = \gamma\circ\beta$, or
an injective contractive linear map $\rho : A\to B$ such that
$\beta\circ\rho = \rho\circ\alpha$. Then $ht(\beta ) > 0$. This conclusion
also holds more generally if $\alpha$ can be obtained from $\beta$ via
a finite chain of intermediary automorphisms intertwined in succession by 
maps of the same form as $\gamma$ or $\rho$.
\end{corollary}

\begin{example}\label{E-operator}
Using Theorem~\ref{T-zero} we can exhibit positive Voiculescu-Brown entropy
in a large class of systems constructed in a operator-theoretic 
fashion as demonstrated by the following examples.
Let $f\in\{ -1,0,1 \}^{\mathbb{Z}}$ be a sequence in which every finite
string of $-1$'s and $1$'s is represented. For each $i=-1,0,1$ we set
$E_i = \{ k\in\mathbb{Z} : f(k) = i \}$. Let $x\in\mathcal{B}
(\ell_2 (E_{-1} \cup E_1 ))$ be the operator obtained by specifying by
$$ x\xi_k = f(k) \xi_k $$
on the set $\{ \xi_k : k\in E_{-1} \cup E_1 \}$ of standard basis elements.
Let $y$ be any self-adjoint operator in 
$\mathcal{B} (\ell_2 (E_0 ))$ of norm at most $1$
and set $a = (x,y) \in\mathcal{B}(\ell_2 (E_{-1} \cup E_1 )) \oplus
\mathcal{B} (\ell_2 (E_0 ))\subset\mathcal{B} (\ell_2 (\mathbb{Z} ))$. Let 
$u$ be the shift $u\xi_k = \xi_{k+1}$ on $\mathcal{B} (\ell_2 (\mathbb{Z} ))$ 
with respect to the canonical basis $\{ \xi_k : k\in\mathbb{Z} \}$, and let
$A\subset\mathcal{B} (\ell_2 (\mathbb{Z} ))$ be the $C^*$-algebra generated 
by $\{ u^n a u^{-n} \}_{n\in\mathbb{Z}}$. By restricting $\Ad\, u$ to $A$
we obtain an automorphism $\alpha$ of $A$. By our assumption on $f$, 
for every $g\in\{ -1 , 1 \}^{\{ 0,\dots , n-1 \}}$ we can find a 
$j\in\mathbb{Z}$ such that for each $k=0, \dots , n-1$ we have 
$a \xi_{j-k} = g(k) \xi_{j-k}$ and hence
$$ \alpha^k (a)\xi_j = g(k) \xi_j . $$
As a consequence, for each $n\in\mathbb{N}$ the 
real linear map which sends the $k$th standard basis element of 
$\ell^n_1$ over the real scalars to $\alpha^k (a)$ for each $k=0, \dots ,
n-1$ is an isometry, and the complexification of this map is an isomorphism
of norm at most $2$ with inverse of norm at most $2$. 
Lemma~\ref{L-basis} then yields
$h_{\topol}(T_\alpha ) = \infty$, and so it follows from Theorem~\ref{T-zero}
that $ht(\alpha ) > 0$ whenever $A$ is exact (this can in fact also be 
deduced directly from Lemma~\ref{L-ell1}---see Remark~\ref{R-cpal}). 
In the case that $a$ is a diagonal operator with respect to the canonical 
basis of $\ell_2 (\mathbb{Z} )$, the $C^*$-algebra $A$ is commutative and the 
topological entropy of the induced homeomorphism of the pure state space
coincides with $ht(\alpha )$ by Proposition 4.8 of \cite{Voi} and hence is 
positive (as can also be seen from Theorem A of \cite{GW}).
The main point of these examples is to demonstrate that positive
Voiculescu-Brown entropy 
can be established in many systems without having either to relate 
the given system to a topological dynamical system which is known a priori to
have positive topological entropy or to rely on measure-theoretic dynamical
invariants like CNT or Sauvageot-Thouvenot entropy. It is sufficient, for 
example, that the eigenspaces of the iterates of an operator of norm $1$ 
corresponding to the respective eigenvalues $\pm 1$ are sufficiently mixed 
along their intersection.
\end{example}

Using Theorem~\ref{T-zero} we can also show that, for the shift on the reduced
crossed product $C^*_r (\mathbb{F}_\infty )$ of the free group on countably
many generators, the topological entropy on the state space is zero
(cf.\ Proposition~\ref{P-full}):

\begin{proposition}\label{P-red}
With $\alpha$ the shift on $C^*_r (\mathbb{F}_\infty )$ we have
$$ h_{\topol}(T_\alpha ) = 0 . $$
\end{proposition}

\begin{proof}
By \cite{Dyk} or \cite{BC} the Voiculescu-Brown entropy of $\alpha$ is 
zero, and so we can apply Theorem~\ref{T-zero}.
\end{proof}

We don't know whether or not the converse of Theorem~\ref{T-zero} holds.
As a test case for this problem we might consider the class of
automorphisms of rotation $C^*$-algebras arising from a
matrix in $SL(2,\mathbb{Z})$ with eigenvalues off the unit circle 
\cite{Wat,Bre}. By \cite{NTA} such a noncommutative $2$-toral automorphism
has positive Voiculescu-Brown entropy, but the arguments in \cite{NTA}
using opposite maps and tensor products give no clue about the value of 
topological entropy on the state space, which we have been unable to
determine in the case of irrational rotation parameters.
Note that for rational rotation parameters
the topological entropy on the state space is infinite since we obtain 
the corresponding commutative $2$-toral automorphism as a subsystem 
by restricting to the centre of the $C^*$-algebra.

We conclude this section with three remarks. 

\begin{remark} 
To show in the proof of Proposition 4.8 in
\cite{Voi} that the Voiculescu(-Brown) entropy dominates the
topological entropy on the pure state space in the separable commutative
setting, Voiculescu applies the classical variational principle along with
several properties of the Connes-Narnhofer-Thirring entropy. It has been a
problem to find a proof of this inequality that does not involve
measure-theoretic entropies. In this regard Theorem~\ref{T-zero} at least 
gives a geometric picture of why positive topological entropy on the
pure state space yields positive Voiculescu-Brown entropy at the 
$C^*$-algebra level. 
\end{remark}

\begin{remark}
As a corollary to Theorem~\ref{T-zero} we recover the result of Glasner and 
Weiss asserting that if a homeomorphism of a compact metric space 
has zero topological entropy then the induced homeomorphism on the space
of probability measures also has zero topological entropy \cite{GW}.
To obtain this corollary we merely need the fact that the topological
entropy of a homeomorphism dominates the Voiculescu-Brown entropy of
the induced $C^*$-algebra automorphism, and this can be established by
a straightforward partition of unity argument (see the proof of
Proposition 4.8 in \cite{Voi}). Now in the geometric
approach of \cite{GW} the construction of the key Banach space map
is most easily managed in the zero-dimensional situation, and indeed
an auxiliary reduction result is invoked to 
handle the general case. Thus by adopting a $C^*$-algebraic viewpoint
we have obtained a functional-analytically more streamlined 
geometric proof of Glasner and Weiss's result.
\end{remark}

\begin{remark}\label{R-cpal}
Lemma~\ref{L-ell1} also yields a means for obtaining lower bounds
for Voiculescu-Brown entropy directly at the completely 
positive approximation level, and is particularly useful when dealing with
dynamical extensions. To illustrate, let $\Omega = \{ x_1 , \dots , x_n \}$ 
be a subset of the unit ball of 
an exact $C^*$-algebra $A$, and suppose that the linear map 
$\Gamma : \ell^n_1 \to\spn\, \Omega$ which sends the $i$th standard basis 
element of $\ell^n_1$ to $x_i$ for each $i=1, \dots , n$ is an isomorphism 
whose inverse is bounded in norm by some $K\geq 1$. Since $\Gamma$ is 
necessarily contractive, it is a $K$-isomorphism, 
i.e., $\| \Gamma \| \| \Gamma^{-1} \| \leq K$. Now if $(\phi , \psi , B)
\in\CPA (\pi , \Omega , \delta )$ for some faithful representation
$\pi : A\to\mathcal{B}(\mathcal{H})$ and $0 < \delta < K^{-1}$, then for any
linear combination $\sum c_i x_i$ of the elements of $\Omega$ we have
\begin{align*}
\Big\| \sum c_i x_i \Big\| &\leq \Big\| \pi\Big( \sum c_i x_i \Big)
- (\psi\circ\phi )\Big( \sum c_i x_i \Big) \Big\|
+ \Big\| (\psi\circ\phi )\Big( \sum c_i x_i \Big) \Big\| \\
&\leq \delta \sum |c_i | + \Big\| \phi \Big( \sum c_i x_i \Big) \Big\| \\
&\leq K\delta \Big\| \sum c_i x_i \Big\| + \Big\| \phi \Big( \sum c_i x_i 
\Big) \Big\| 
\end{align*}
so that $\big\| \phi \big( \sum c_i x_i \big) \big\| \geq (1-K\delta )
\big\| \sum c_i x_i \big\|$,
and since $\phi$ is contractive it follows that $\phi | {}_{\spn\, \Omega}$ is 
a $(1-K\delta )^{-1}$-isomorphism onto its image. Hence the Banach-Mazur
distance between $\ell^n_1$ and the image of $\spn\, \Omega$ under $\phi$
is at most $K(1-K\delta )^{-1}$, and so by Lemma~\ref{L-ell1} we conclude that
$$ \log\rcp (\Omega , \delta ) \geq n a^{-1} K^{-2} (1-K\delta )^2 $$
where $a > 0$ is a universal constant.
Thus if $\alpha$ is an automorphism of $A$ and $x$ is an element of $A$
such that there exists a $K\geq 1$ such that for every $n\in\mathbb{N}$ the 
linear map which sends the $i$th standard basis element of $\ell^n_1$ to 
$\alpha^{i-1} (x)$ for each $i=1, \dots , n$ is an isomorphism whose inverse 
has norm at most $K$, then
$$ ht(\alpha ) \geq\sup_{\delta > 0} ht(\alpha , \{ x \} , \delta )\geq
\sup_{\delta\in (0,K^{-1})} a^{-1} K^{-2} (1-K\delta )^2 = a^{-1}K^{-2} . $$
This lower bound for entropy also applies to any automorphism $\beta$ of an 
exact $C^*$-algebra $D$ such that there exists a surjective $^*$-homomorphism 
$\gamma : D\to A$ with $\gamma\circ\beta = \alpha\circ\gamma$, for
in such a case we can lift $x$ under $\gamma$ to a 
element $y\in D$ of the same norm, and the linear map 
sending the $i$th standard basis element of $\ell^n_1$ to 
$\beta^i (y)$ for each $i=1, \dots , n$ is an isomorphism with inverse of 
norm at most $K$, as is easily checked. 
More generally, we also obtain 
$ht(\beta ) > 0$ for any automorphism $\beta$ of a separable 
exact $C^*$-algebra such that $\alpha$ can obtained from $\beta$
via a finite chain of intermediary automorphisms intertwined in
succession by contractive linear surjections or 
linear isometries in the reverse direction (cf.\ Corollary~\ref{C-pos}). 
\end{remark}

\section{A description of the topological Pinsker 
algebra in terms of local Voiculescu-Brown entropy}\label{S-Pinsker}

Let $T:X\to X$ be a homeomorphism of a compact metric space and 
$\alpha_T$ the automorphism of $C(X)$ given by $\alpha_T (f)
= f\circ T$ for all $f\in C(X)$. We recall from \cite{Disj} that a
pair $(x,y)\in X\times X\setminus\Delta$ (with $\Delta$ denoting the diagonal)
is called an {\it entropy pair} if $h_{\topol} (T,\mathcal{U}) > 0$
for every two-element open cover $\mathcal{U} = \{ U,V \}$ with
$x\in\interior (X\setminus U)$ and $x\in\interior (X\setminus V)$. 
We denote by $E_X$ the set of entropy pairs in $X\times X$. 
The {\it topological Pinsker factor} is defined as the quotient system 
arising from the closed $T$-invariant equivalence relation on $X$ generated 
by the collection of entropy pairs \cite{BL} (here we have adopted the 
terminology of \cite{Gl}). This translates at the 
$C^*$-algebra level as the $\alpha_T$-invariant $C^*$-subalgebra $P_{X,T}$
of $C(X)$ consisting of all $f\in C(X)$ satisfying $f(x) = f(y)$ for every
entropy pair $(x,y)$. Note that $P_{X,T}$ is indeed $\alpha_T$-invariant 
because $E_X$ is invariant under $T\times T$ by Proposition 3 of \cite{Disj}. 
We refer to $P_{X,T}$ as the
{\it topological Pinsker algebra}. It is an analogue of the Pinsker 
$\sigma$-algebra in ergodic theory (see \cite{Wal}).

The main goal of this section is to apply the argument of the proof 
of Theorem~\ref{T-zero} to show that $P_{X,T}$ is equal to the set of all 
$f\in C(X)$ such that the local Voiculescu-Brown entropy of $\alpha_T$ with
respect to the singleton $\{ f \}$ is zero. Thus by viewing 
the dynamics at the function level we are able to obtain a simple 
description of the topological Pinsker factor/algebra that avoids entropy
pairs and the awkward fact 
that the set $E_X \cup\Delta$ does not always form an equivalence relation
(see \cite{TEE}). As shown in \cite{NTA}, our functional-analytic 
description of $P_{X,T}$ can be applied to obtain some information
concerning the positivity of local Voiculescu-Brown entropy with
respect to products of canonical unitaries for certain
noncommutative toral automorphisms.

For economy, in this section we will simply write 
$ht(T,f)$ instead of $ht(\alpha_T , \{ f \} )$ 
(as it appears in the definition of Voiculescu-Brown entropy) for any function 
$f\in C(X)$. Also, given a function $f\in C(X)$ we define the
pseudo-metric $d_f$ on $X$ by
$$ d_f (x,y) = | f(x) - f(y) | $$
for all $x,y\in X$.
We furthermore need to extend the metric space formulation
of topological entropy to pseudo-metrics. Thus for any pseudo-metric
$d$ on $X$ we set
$$ h_d (T) = \sup_{\varepsilon >0} \limsup_{n\to\infty}\frac1n
\log\sep_n (T,\varepsilon ) $$
where $\sep_n (T,\varepsilon )$ is the largest cardinality of an 
$(n,\varepsilon )$-separated set, with the latter defined in the same way as 
for metrics (see Section~\ref{S-top}).

The following lemma is a local version of Theorem~\ref{T-zero} at the
level of a single element in the $C^*$-algebra.

\begin{lemma}\label{L-locmetric}
Let $f\in C(X)$. Then $h_{d_f}(T) > 0$ implies $ht(T,f) > 0$.
\end{lemma}

\begin{proof}
Notice that, in the proof of Theorem~\ref{T-zero}, if the compact set
$K$ is taken to be finite, then the argument 
still shows that if we define the pseudo-metric
$$ d(\sigma , \omega ) = \max_{x\in K} | \sigma (x) - \omega (x) | $$
on the unit ball of the dual
and take $h_d (S_\alpha ) > 0$ as our hypothesis, then
$ht(\alpha , K) > 0$. Thus in our present context we can take $K$ in the 
proof of Theorem~\ref{T-zero} to be the singleton
$\{ f \}$ to obtain the desired conclusion.
\end{proof}

\begin{lemma}\label{L-epair}
If $(x,y)$ is an entropy pair then 
$ht(T,f) > 0$ for every $f\in C(X)$ with $f(x) \neq f(y)$.
\end{lemma}

\begin{proof}
Suppose $(x,y)$ is an entropy pair and $f$ is a function in $C(X)$ with
$f(x) \neq f(y)$. Set $\delta = | f(x) - f(y) |/3$ and define the two 
open sets
\begin{align*}
U &= \{ z\in X : | f(x) - f(z) | > \delta \} , \\
V &= \{ z\in X : | f(x) - f(z) | < 2\delta \} .
\end{align*}
Then $\mathcal{U} = \{ U,V \}$ is an open cover with $x\in\interior 
(X\setminus U)$ and $y\in\interior (X\setminus V)$, and so $h_{\topol} 
(T,\mathcal{U}) > 0$ by virtue of the fact that $(x,y)$ is an entropy pair.
Now if $n\in\mathbb{N}$ and $\mathcal{V}$ is a subcover of 
$\mathcal{U}\vee T^{-1} \mathcal{U} \vee\cdots\vee T^{-(n-1)} \mathcal{U}$
of smallest cardinality then in each element of $\mathcal{V}$ we can
choose a point which is not contained in any other element of $\mathcal{V}$,
for otherwise $\mathcal{V}$ would not be a minimal subcover. The set $E$
obtained by collecting these points together is $(n,\delta/2)$-separated
relative to the pseudo-metric $d_f$ and has the same cardinality as 
$\mathcal{V}$. Hence
$$ h_{d_f}(T) \geq h_{\topol} (T,\mathcal{U}) > 0 $$
and so $ht(T,f) > 0$ by Lemma~\ref{L-locmetric}, yielding the result. 
\end{proof}

The converse of Lemma~\ref{L-epair} is false. This is a consequence of the
fact that set $E_X \cup\Delta$ is not necessarily transitive as a 
relation (see \cite{TEE}). Indeed if $(x,y)$ and $(y,z)$ are elements of
$E_X \cup\Delta$ such that $(x,z)\not\in E_X \cup\Delta$, and if $f\in C(X)$ 
satisfies $f(x) \neq f(z)$, then we cannot have both $f(x) = f(y)$ and 
$f(y) = f(z)$, whence $ht(T,f) > 0$ by Lemma~\ref{L-epair}.

\begin{theorem}\label{T-Pinsker}
The topological Pinsker algebra $P_{X,T}$ is equal to the set of all 
$f\in C(X)$ such that $ht(T,f) = 0$.
\end{theorem}

\begin{proof}
What we need to prove is that for any $f\in C(X)$ we have $ht(T,f) = 0$
if and only if $f(x) = f(y)$ for all entropy pairs $(x,y)$. The ``only if''
direction follows immediately from Lemma~\ref{L-epair}. For the ``if''
direction, let $f\in C(X)$ and suppose $ht(T,f) > 0$. Let $B$ be the unital 
$\alpha_T$-invariant $C^*$-subalgebra of $C(X)$ generated by 
$\{ \alpha^n_T (f) \}_{n\in\mathbb{Z}}$,
i.e., the closure in $C(X)$ of the set of polynomials in
$\{ \alpha^n_T (f) \}_{n\in\mathbb{Z}}$. For convenience, in the rest of the
proof we will identify points in a compact metric space with pure states on
the corresponding unital $C^*$-algebra. Now if $T_B$ denotes the 
homeomorphism induced by $\alpha_T \big| {}_B$
on the pure state space of $B$, then since $ht(\alpha_T \big| {}_B ) \geq
ht(T,f) > 0$ we have $h_{\topol} (T_B ) > 0$ by \cite[Prop.\ 4.8]{Voi}, 
and so by \cite[Props.\ 1 and 2]{Disj} there are pure
states $\sigma , \omega$ on $B$ such that $(\sigma , \omega )$ is an
entropy pair with respect to $T_B$.
We must then have $\sigma (\alpha^n_T (f)) \neq\omega (\alpha^n_T (f))$
for some $n\in\mathbb{Z}$ since $\sigma$ and $\omega$ are distinct, and
since the set of entropy pairs is $T\times T$-invariant 
\cite[Prop.\ 3]{Disj} we may
assume that $n=0$, i.e., $\sigma (f) \neq\omega (f)$. By 
\cite[Prop.\ 4]{Disj} there are pure states $\sigma'$ and $\omega'$ 
on $C(X)$ extending $\sigma$ and $\omega$, respectively, such that 
$(\sigma' , \omega' )$ forms an entropy pair, and we
have $\sigma' (f) \neq\omega' (f)$, completing the proof.
\end{proof} 
 
The system $(X,T)$ is said to have {\it completely positive entropy} if 
each of its non-trivial factors has positive topological 
entropy \cite{FPTE}.
Since positive entropy systems always have an entropy pair 
\cite[Props.\ 1 and 2]{Disj}, $P_{X,T}$ is equal to the scalars 
(resp.\ $C(X)$) precisely when the system
$(X,T)$ has completely positive entropy (resp.\ zero entropy),
and so we obtain the following corollaries to 
Theorem~\ref{T-Pinsker}.

\begin{corollary}\label{C-cpe}
The system $(X,T)$ has completely positive entropy if and only if 
$ht(T,f) > 0$ for all non-constant functions $f\in C(X)$.
\end{corollary}

\begin{corollary}\label{C-posentr}
We have $h_{\topol} (T) > 0$ if and only if
$ht(T,f) > 0$ for some $f\in C(X)$.
\end{corollary}

\end{document}